\documentclass{siamltex}
\usepackage{lineno,hyperref}
\usepackage{multirow}
\usepackage[english]{babel}
\usepackage{amsmath}
\usepackage{bm}
\usepackage{ulem}
\usepackage[matha,mathx]{mathabx}
\usepackage{mathrsfs}
\usepackage{amssymb}
\usepackage{graphicx}
\usepackage{tikz}
\usepackage{here}
\usepackage{algorithm}
\usepackage{algorithmic}
\usepackage{graphicx}
\usepackage{amsmath}
\usepackage{here}
\usetikzlibrary{matrix,calc}
\usepackage{caption}
\usepackage{enumitem}

\newcommand{\vs}{\vspace{0.2cm}}
\makeatletter
\newcommand*{\rom}[1]{\expandafter\@slowromancap\romannumeral #1@}
\makeatother

\usepackage{tikz}
\usepackage{here}
\usepackage{algorithm}
\usepackage{algorithmic}
\usepackage{graphicx}   
\usepackage{amsmath}
\usepackage{here}
\usetikzlibrary{matrix,calc}
\usepackage{caption}
\DeclareMathAlphabet{\mathpzc}{OT1}{pzc}{m}{it}
\title{Tensor Golub Kahan based on Einstein product}

\author{  A. El Hachimi\footnotemark[2] \thanks{The UM6P Vanguard Center, Mohammed VI Polytechnic University, Green 
City, Morocco.}
\and K. Jbilou\footnotemark[1] \thanks{Universit\'e du Littoral Cote d'Opale, LMPA, 50 rue F. Buisson, 62228 Calais-Cedex, France.}
\and M. Hached \thanks{Universit\'e de Lille, CNRS, UMR 8524 - Laboratoire Paul Painlev\'e, F-59000 Lille, France.}
\and A. Ratnani\footnotemark[1]}

\begin{document}
\maketitle 

\begin{abstract}
The Singular Value Decomposition (SVD) of matrices is a widely used tool in scientific computing. In many applications of machine learning, data analysis, signal and image processing, the large datasets are structured into tensors, for which generalisations of SVD have already been introduced, for various types of tensor-tensor products. In this article, we present innovative methods for approximating this generalisation  of SVD to tensors in the framework of the Einstein tensor product. These singular elements are called singular values and singular tensors respectively. The proposed method uses the tensor Lanczos bidiagonalization applied to the Einstein product. In most applications, as in the matrix case, the extremal singular values are of special interest. To enhance the approximation of the largest or the smallest singular triplets (singular values and left and right singular tensors), a  restarted method based on Ritz augmentation is proposed. Numerical results are proposed to illustrate the effectiveness of the presented method. In addition, applications to video compression and facial recognition are presented.
\end{abstract}

\section{introduction}

One of the ubiquitous and challenging problems in mathematics is the extraction of the most relevant information out of large datasets, especially in the domains of data mining, machine learning and deep learning. In the context where data is represented by matrices, the determination of the main features  is closely related to the Singular Value Decomposition (SVD) of some matrix, like the covariance matrix for Principal Component Analysis (PCA). More precisely, the main features of the data are related to the  largest singular values of a matrix and their associated left and right singular vectors, summarised under the name of singular triplets. Due the generally large size of the datasets, the computation of the SVD is sometimes unrealistic in terms of computational time. This problem has already been addressed in the matrix case by Golub and al. in \cite{golub}, based on a Lanczos bidiagonalization (\textbf{LB}) algorithm. This approach, which was later refined by Baglama and Reichel in \cite{baglama2005augmented}, aims to avoid the computation of the SVD of a potentially very large matrix, approximating just a chosen number of the largest (or smallest) singular triplets.

\noindent In numerous real world applications, due to the large quantity of information and variables, datasets are often structured as multidimensional arrays, the so-called \textbf{\textit{tensors}}.  A tensor is a multidimensional array of numbers, and can be seen as a generalisation of matrices and vectors. An $N$-order tensor (or $N$-mode tensor) is a $N$ dimensional array (a vector for $N=1$ and a matrix for $N=2$). Various tensor-tensor products have been defined, generalising the matrix-matrix product, as the t-product, and the c-product that have been introduced by Kilmer et all. \cite{kilmer1,kilmer2013third} and the most notorious Einstein product \cite{Einstein}.  The t and c-products are defined via matrix-matrix products in the Fourier and Cosine domain, after a Fast Fourier transform (FFT) or Discrete Cosine transform (DCT). These products were first designed only for third-order tensors, although Martin et al. and  Bentbib and al. have generalised them for higher order tensor in \cite{Martin, elha4}. Note that in their article \cite{elha4}, the authors also gave a more  general tensor-tensor product using any invertible linear transform, the $\mathcal{L}$-product.

In the recent years, tensor-based numerical methods have received a great deal of attention for their ability to address very topical problems. Amongst many other applications, completion problems, which consist of approximating some missing information in large datasets (recommendation systems for customer, restoration of damaged images etc.), see \cite{elha1,elha3} for more details. We can also mention the Tensor Robust Principal Component Analysis (RPCA), which objective is to remove the noise from corrupted data \cite{elha2}, image restoration \cite{elguide2,BKS2} and numerical meshless methods for the resolution of 3-d PDE's via Radial Basis Functions (RBF) \cite{elguide1}.

\noindent The generalisation of the SVD for tensors has been introduced in various ways, depending on the chosen tensor-tensor product. The t-svd, c-svd, and $\mathcal{L}$-svd in \cite{kilmer1,kernfeld2015tensor,elha3} are defined via the t-product, the c-product and the $\mathcal{L}$-product, respectively, and have been originally designed for 3-mode tensors. In this paper, we will focus on the SVD associated to the more general Einstein product \cite{Einstein2013,Einstein2026}. 

\noindent The approximation of the largest singular triplets for  tensors was first introduced in \cite{Hached} for 3-mode tensors, using the c-product. This approach was an adaptation to tensors of the matrix Lanczos bidiagonalization. In order to improve the accuracy of this approximation, in \cite{elha4}, the authors generalized the work done by Baglama and Reichel \cite{baglama2005augmented} under the t-product, by restarted strategies using Ritz and Harmonic Ritz tensors. This approach gave very satisfactory results in terms of accuracy and were illustrated by applications to facial recognition, image classification and data compression \cite{elha4,Hached}.

\noindent This paper is outlined as follows: In Section \ref{sec 2}, a reminder of the Einstein product and some properties are presented, Section \ref{sec 3} is devoted to the approximation of the tensor singular triplets can  by using the tensor Lanczos bidiagonalization associated with the Einstein product, along with a restart strategy by Ritz tensors. Section \ref{sec PCA} presents a new Einstein product based PCA technique, and some numerical experiments are presented in Section \ref{sec 5}.
\section{Definitions and notations}
\label{sec 2}
\noindent In this section, we review some notations, definitions, and properties linked with tensors. The notations used throughout this paper are the same as those employed by Kolda and Bader in \cite{kolda2009tensor}. A real $N$th-order tensor $\mathcal{A}$ is an $N$-dimensional array of real numbers and we denote $\mathbb{R}^{I_1\times \ldots \times I_N}$ the real vector space of $N$th-order real tensors. The $(i_1,\ldots, i_N)^{\text{th}}$ entry of the tensor $\mathcal{A}$ is denoted by $\mathcal{A}_{i_1\ldots i_N}$. General definitions about foldings, unfoldings, fibers and slices can be found in \cite{kolda2009tensor}.

\vs
\noindent First, we describe the n-mode product which defines a multiplication of a tensor by a matrix or a vector.
\vs
\begin{definition}\cite{kolda2009tensor}
Let $\mathcal{A}\in \mathbb{R}^{I_1\times \ldots \times I_N}$ and $U\in \mathbb{R}^{J\times I_n}$ be an $N$-order tensor and a matrix, respectively. The n-mode product between $\mathcal{A}$ and $U$ produces a tensor of size $I_1\times \ldots \times I_{n-1}\times J\times I_{n+1}\times \ldots \times I_N$, which entries are given by
\[
\left(\mathcal{A}\times_n U\right)_{i_1\ldots i_{n-1}ji_{n+1}\ldots i_N}=\sum_{i_n=1}^{I_n}\mathcal{A}_{i_1\ldots i_N} U_{j i_n}.
\]
\end{definition}

\noindent We recall the following identities :\\

\vs
\noindent For all tensor $\mathcal{A}\in \mathbb{R}^{I_1\times \ldots \times I_N}$ and for all matrices $U\in \mathbb{R}^{J\times I_n}$ and $V\in \mathbb{R}^{K\times I_m}$, we have
\[
\mathcal{A}\times_n U \times_m V =\mathcal{A}\times_m V \times_n U.
\]
\noindent Moreover, if  $U\in \mathbb{R}^{J\times I_n}$ and $V\in \mathbb{R}^{I_n \times J}$, then, we have 
\[
\mathcal{A}\times_n U\times V=\mathcal{A}\times_n\left(VU\right).
\]
\noindent We now give the definition of the n-mode matricization of a tensor. 

\vs 

\begin{definition}\cite{kolda2009tensor}
Let $\mathcal{A}\in \mathbb{R}^{I_1\times \ldots  \times I_N}$  be an $N$-order tensor. The n-mode matricization of $\mathcal{A}$ is denoted by $\mathcal{A}_{(n)}$ belongs to $\mathbb{R}^{I_n \times \prod_{i=1}^{N}I_i}$. It maps the $(i_1,\ldots,i_N)$-th element fo $\mathcal{A}$ to the $(i_n,j)$-th element of $\mathcal{A}_{(n)}$ one, such that
\[
j=1+\sum_{k=1,k\neq n }^{N}(i_k-1)J_k, \quad \text{with}\quad J_k=\prod_{m=1,m\neq m}^{k-1}I_m.
\]
\end{definition}

\noindent The following identity gives an interpretation of the n-mode product in terms of a matrix product:

\vs

\noindent Let $\mathcal{A}\in \mathbb{R}^{I_1\times \ldots \times I_N}$ and $B\in \mathbb{R}^{J\times I_n}$, then we have
\begin{equation}\label{n-mode}
\mathcal{C}=\mathcal{A}\times_n B \Longleftrightarrow \mathcal{C}_{(n)}=B\mathcal{A}_{(n)},
\end{equation}
where $\mathcal{C}_{(n)}$ and $\mathcal{A}_{(n)}$ are respectively the n-mode matricization of $\mathcal{A}$ and  $\mathcal{C}$.

\vs

Next, we recall the Einstein product \cite{Einstein,Einstein2013,Einstein2026}. 

\vs

\begin{definition}\cite{Einstein2013}
Let $\mathcal{A}\in \mathbb{R}^{I_1\times \ldots\times I_L \times K_1 \times \dots \times K_N}$ and $\mathcal{B}\in \mathbb{R}^{K_1\times \ldots \times K_N\times J_1\times \ldots \times J_N}$. The Einstein product between $\mathcal{A}$ and $\mathcal{B}$ is the tensor of size $I_1\times \ldots \times I_L \times J_1\times \ldots \times J_N$ whose elements are defined by 
\[
\left(\mathcal{A}*_N\mathcal{B}\right)_{i_1\ldots i_L j_1\ldots j_M}=\sum_{k_1,\ldots,k_N=1}^{K_1,\ldots, K_N}\mathcal{A}_{i_1\ldots i_L k_1\ldots k_N}\mathcal{B}_{k_1\ldots k_N j_1\ldots j_M}.
\]
\end{definition}

\vs

\noindent The following definitions will be needed in the sequel.

\vs

\begin{definition}
\begin{itemize}
	\item \textbf{Transpose tensor:} For a given tensor $\mathcal{A}\in \mathbb{R}^{I_1\times \ldots \times I_N\times J_1\times \ldots \times J_M}$ the transpose tensor $\mathcal{B}$ of $\mathcal{A}$ denoted by $\mathcal{A}^T$ is the tensor of size $J_1\times \ldots \times J_M \times I_1\times \ldots \times I_N$ whose elements are $\mathcal{B}_{j_1\dots j_M i_1\dots i_N}=\mathcal{A}_{i_1\dots i_N j_1\dots j_M}$.
	\item \textbf{Diagonal tensor:} A tensor $\mathcal{D}\in \mathbb{R}^{I_1\times \ldots \times I_N \times J_1\times \ldots \times J_M}$ is a diagonal tensor if\\ $\mathcal{D}_{i_1\dots i_N j_1\dots j_M}= 0$ for all $i_k\neq j_k$ with $k\in\{1,\ldots,\min(N,M)\}$.
	\item \textbf{Identity tensor:} $\mathcal{I}_N\in \mathbb{R}^{I_1\times \ldots \times I_N \times I_1 \times \ldots \times I_N}$ is called an identity tensor if it is diagonal and all its diagonal entries are ones.
	\item \textbf{Invertible tensor:} A tensor $\mathcal{A}\in \mathbb{R}^{I_1\times \ldots I_N \times I_1\times \ldots \times I_N}$ is invertible if and only if there exists a tensor $\mathcal{X}\in \mathbb{R}^{I_1\times \ldots I_N \times I_1\times \ldots \times I_N}$ such that 
	\[
	\mathcal{A}*_N\mathcal{X}=\mathcal{X}*_N \mathcal{A}=\mathcal{I}_N.
	\]
	If such a tensor exists, it is unique. It is called the inverse of $\mathcal{A}$ and is denoted by $\mathcal{A}^{-1}$.
\end{itemize}
\end{definition}

\vs

\noindent Let $\mathcal{A}, \mathcal{B}$ two tensors in $\mathbb{R}^{I_1\times \ldots \times I_N \times J_1\times \ldots \times J_M}$. The inner product between $\mathcal{A}$ and $\mathcal{B}$ is defined by 
\[
\langle \mathcal{A},\mathcal{B}\rangle = {\tt tr}\left(\mathcal{A}^T *_N \mathcal{B}\right),
\]
where, the trace of $\mathcal{A}\in \mathbb{R}^{I_1\times \ldots \times I_N \times I_1\times \ldots \times I_N}$ is defined as
\[
{\tt tr}\left(\mathcal{A}\right)=\sum_{i_1,\ldots,i_N=1}^{I_1,\ldots,I_N}\mathcal{A}_{i_1\ldots i_N i_1\dots i_N}.
\]
If we have $\langle \mathcal{A},\mathcal{B}\rangle =0$, then we say that $\mathcal{A}$ and $\mathcal{B}$ are orthogonal.\\
The corresponding norm is the tensor Frobenius norm given by
\[
\left\Vert \mathcal{A}\right\Vert_F=\sqrt{{\tt tr}\left(\mathcal{A}^T *_N \mathcal{A}\right)}.
\]

\begin{proposition}\cite{Einstein2026}
Let $\mathcal{A}\in \mathbb{R}^{I_1\times \ldots \times I_L \times K_1\times \ldots \times K_N}$ and $\mathcal{B}\in \mathbb{R}^{K_1\times \ldots \times K_N\times J_1\times \ldots \times J_M}$. Then 
\[
\left(\mathcal{A}*_N \mathcal{B}\right)^T=\mathcal{B}^T *_N \mathcal{A}^T.
\]
\end{proposition}
\noindent In the sequel, for all considered tensor $\mathcal{A}\in \mathbb{R}^{I_1\times \ldots \times I_N \times J_1 \times \ldots \times J_M}$, we will suppose that $M\leq N$.

\vs

\noindent The generalization of the svd under the Einstein product is defined as follows :

\vs

\begin{theorem}\cite{Einstein2026}
\noindent Let us consider a tensor $\mathcal{A}\in \mathbb{R}^{I_1\times \ldots \times I_N \times J_1 \times \ldots \times J_M}$. The singular value decomposition (SVD) of $\mathcal{A}$ is given by a decomposition
\[
\mathcal{A}=\mathcal{U}*_N \mathcal{S}*_M \mathcal{V}^T,
\]
where $\mathcal{U}\in \mathbb{R}^{I_1\times \ldots \times I_N\times I_1\times \ldots \times I_N}$, $\mathcal{V}\in \mathbb{R}^{J_1\times \ldots \times J_M\times J_1\times \ldots \times J_M}$, and $\mathcal{S}\in \mathbb{R}^{I_1\times \ldots \times I_N\times J_1\times \ldots \times J_M}$ is a diagonal tensor satisfying:
\begin{itemize}
	\item If $M=N$
	\[
	\mathcal{S}_{i_1\ldots i_N j_1\ldots j_N}=\left\lbrace\begin{array}{ll}
		\mu_{i_1 \ldots i_N}  \text{~if~~}(i_1,\ldots,i_N)=(j_1,\ldots,j_N),\\
		0 \text{~if~~} (i_1,\ldots,i_N)\neq (j_1,\ldots,j_N).
	\end{array}\right.
	\]
	\item If $M<N$
	\[
	\mathcal{S}_{i_1\ldots i_N  j_1\ldots j_M}=\left\lbrace\begin{array}{ll}
		\mu_{i_1 \ldots i_M} \text{~if~~} (i_1,\ldots,i_M)=(j_1,\ldots,j_M), \; (i_{M+1},\ldots, i_N)=(1,\ldots, 1),\\
		0  \quad \text{otherwise}.
	\end{array}\right.
	\]
\end{itemize}
\noindent where the positive real numbers $\mu_{i_1 \ldots i_M}$ are called the singular values of $\mathcal{A}$.
\end{theorem}

\vs

\noindent Assume that we have $\mathcal{A}\in \mathbb{R}^{I_1\times \ldots \times I_N\times J_1\times \ldots \times J_M}$ and let $\widetilde{I}_i=\min(I_i,J_i)$ for $i=1,\ldots,M$. Let us consider the tensor svd given in the above theorem. Therefore, we can write :
\[
\mathcal{A}=\sum_{i_1,\ldots, i_M=1}^{\widetilde{I}_1,\ldots,\widetilde{I}_M} \mu_{i_1 \ldots i_M}\mathcal{U}_{i_1\ldots i_M}*_M \left(\mathcal{V}_{i_1\ldots i_M}\right)^T,
\]
\noindent where $\mathcal{U}_{i_1 \ldots i_M} \in \mathbb{R}^{I_1\times \ldots \times I_N}$ and $\mathcal{V}_{i_1 \ldots i_M} \in \mathbb{R}^{J_1\times \ldots \times J_N}$ are called the left and the right singular tensors respectively. Note that we have :

\vs

\begin{itemize}
\item If $M=N$, $\mathcal{U}_{i_1\ldots i_M}=\mathcal{U}(:,\ldots,:, i_1,\ldots, i_M)$,
\item If $M< N$, $\mathcal{U}_{i_1\ldots i_M}=\mathcal{U}(:,\ldots,:, i_1,\ldots, i_M, 1,\ldots , 1)$,
\item $\mathcal{V}_{i_1 \ldots i_M}=\mathcal{V}(:,\ldots,:, i_1,\ldots, i_M)$.
\end{itemize} 

\vs

\noindent The singular values $\mu_{i_1\ldots I_N}$ are ordered decreasingly, \textit{i.e.},
\[
\mu_{1,1,\ldots,1} \geq \mu_{2,1,\ldots,1} \geq \dots \geq  \mu_{\widetilde{I}_1,1,\ldots,1} \geq  \ldots  \geq \mu_{\widetilde{I}_1,\widetilde{I}_2,\ldots, \widetilde{I}_N}.
\]
Notice that $\mathcal{A}$ admits $\widetilde{I}_1\widetilde{I}_2\ldots \widetilde{I}_N$ singular values and left and  right singular tensors.\\
In the sequel, for $k\in \mathbb{N}^*$, $s_k$,  $\mathcal{U}_k$ and $\mathcal{V}_k$ will denote the $k^{th}$  singular value and the left and right singular tensors respectively and  $\left\{s_k, \mathcal{U}_k, \mathcal{V}_k\right\}$ will be called the $k^{th}$ singular triplet of $\mathcal{A}$.
\section{The large tensor svd based on tensor Krylov subspace}
\label{sec 3}

The tensors involved in many problems are large, making the computation of the SVD computationally very challenging. Nevertheless, in almost all applications, only the largest singular triplets are required. The well-known Golub-Kahan method has been used to address this issue in the matrix case and very useful refinements have been proposed by Reichel and Baglama \cite{baglama2005augmented}. More recently, in \cite{elha4,Hached}, the authors proposed a generalization of those refinements to the t/c or $L$-products for 3th-order tensors, with applications to PCA, image recognition and classification problems. In this section, we propose an approximation method of singular triplets of large tensors under the Einstein product which is based on the tensor Lanczos bidiagonalization.
\subsection{Einstein tensor Lanczos bidiagonalization}
Let $\mathcal{A}\in \mathbb{R}^{I_1\times \ldots \times I_N \times J_1\times \ldots \times J_M}$ be a large tensor, $\mathcal{P}\in \mathbb{R}^{J_1\times \ldots \times J_N}$ and $m\in \mathbb{N}^*$such as $m\leq \min(I_1\ldots I_N, J_1 \ldots J_M)$. An $m$ tensor Lanczos diagonalization (Golub Kahan) process consists in constructing two orthogonal subspaces from the following Krylov subspaces
\begin{eqnarray*}
\textbf{K}_m\left(\mathcal{A}^T*_N \mathcal{A}, \mathcal{P}\right)={\tt span}\left\{\mathcal{P}, \left(\mathcal{A}^T*_N \mathcal{A}\right)*_M \mathcal{P}, \ldots, \left(\mathcal{A}^T*_N \mathcal{A}\right)^{m-1}*_M \mathcal{P}\right\},\\
\textbf{K}_m\left(\mathcal{A}*_M \mathcal{A}^T, \mathcal{Q}\right)={\tt span}\left\{\mathcal{Q}, \left(\mathcal{A}*_M \mathcal{A}^T\right)*_N \mathcal{Q}, \ldots, \left(\mathcal{A}*_M \mathcal{A}^T\right)^{m-1}*_N \mathcal{Q}\right\},
\end{eqnarray*}
with $\left\Vert \mathcal{P}\right\Vert_F=1$, and $\mathcal{Q}=\mathcal{A}*_N \mathcal{P}$.
\vs
The Einstein-product-based tensor Lanczos bidiagonalization algorithm is summarized as follows :

\begin{algorithm}[H]
\caption{Einstein tensor Lanczos bidiagonalization (ELB)}
\label{alg LB}
\textbf{Input:} $\mathcal{A}\in \mathbb{R}^{I_1\times \ldots \times I_N \times J_1\times \ldots \times J_M}$, $\mathcal{P}_1\in \mathbb{K}^{J_1\times \ldots \times J_M}$ unitary ($\left\Vert \mathcal{P}_1\right\Vert_F=1$) and $0<m\leq\min(I_1\ldots I_N, J_1 \ldots J_M)$.\\
\textbf{Output:} $\mathbb{P}_m \in \mathbb{R}^{J_1\times \ldots \times J_M\times m}$,  $\mathbb{Q}_m \in \mathbb{R}^{I_1\times \ldots \times I_N\times m}$, $B_m\in \mathbb{R}^{m\times m}$, $\mathcal{R}_m\in \mathbb{R}^{J_1\times \ldots \times J_M}$ .

\begin{algorithmic}[1]
	\STATE $\mathcal{Q}_1=\mathcal{A}*_M\mathcal{P}_1$, $\alpha_1=\left\Vert \mathcal{Q}_1\right\Vert_F$, $\mathcal{Q}_1=\dfrac{\mathcal{Q}_1}{\alpha_1}$.
	\FOR{$j=1,\ldots,m$}
	\STATE $\mathcal{R}_j=\mathcal{A}^T*_N \mathcal{Q}_j - \alpha_j \mathcal{P}_j$.
	\IF{$j<m$}
	\STATE $\beta_j=\left\Vert \mathcal{R}_j\right\Vert_F$.
	\IF{$\beta_j=0$}
	\STATE Stop
	\ELSE
	\STATE $\mathcal{P}_{j+1}=\dfrac{\mathcal{R}_j}{\beta_j}$.
	\ENDIF
	\STATE $\mathcal{Q}_{j+1}=\mathcal{A}*_M \mathcal{P}_{j+1}-\beta_j \mathcal{Q}_j$.
	\STATE $\alpha_{j+1}=\left\Vert \mathcal{Q}_{j+1}\right\Vert_F$.
	\IF{$\alpha_{j+1}=0$}
	\STATE Stop
	\ELSE
	\STATE $\mathcal{Q}_{j+1}=\dfrac{\mathcal{Q}_{j+1}}{\alpha_{j+1}}$.
	\ENDIF
	\ENDIF
	\ENDFOR
\end{algorithmic}
\end{algorithm}

\noindent The above algorithm returns the tensors $\mathbb{P}_m$ and $\mathbb{Q}_m$ which are built on the tensors $\mathcal{P}_1,\ldots, \mathcal{P}_m$ and $\mathcal{Q}_1,\ldots , \mathcal{Q}_m$, respectively, \textit{i.e.}, $\mathbb{P}_m=\left[ \mathcal{P}_1,\ldots ,\mathcal{P}_m \right]\in \mathbb{R}^{J_1\times \ldots \times J_M\times m}$ and $\mathbb{Q}_m=\left[   \mathcal{Q}_1,\ldots ,\mathcal{Q}_m \right]\in \mathbb{R}^{I_1\times \ldots \times I_N\times m}$. Using the Matlab notations, we can see that $\mathbb{P}_m(:,\ldots,:, i)=\mathcal{P}_i$ and $\mathbb{Q}_m(:,\ldots,:, i)=\mathcal{Q}_i$. Algorithm \ref{alg LB} returns the bidiagonal matrix $B_m \in \mathbb{R}^{m\times m}$ defined by 
\[
B_m=\begin{pmatrix}
\alpha_1 & \beta_1 &  &&\\
0 & \alpha_2 & \beta_2 & &\\
& & \ddots & \ddots  & \\
& & & \alpha_{m-1} & \beta_{m-1}\\
& & & & \alpha_m
\end{pmatrix}.
\]
It should be mentioned that $\mathcal{P}_1,\ldots, \mathcal{P}_m$ are mutually orthonormal, \textit{i.e.}, $\langle \mathcal{P}_i,\mathcal{P}_j \rangle=\delta_{i,j}$, as well as $\mathcal{Q}_1,\ldots, \mathcal{Q}_m$. This can also be expressed  as $\mathbb{P}_m^T *_M \mathbb{P}_m= I_m$ and $\mathbb{Q}_m^T *_N \mathbb{Q}_m= I_m$.\\

\vs

\noindent The relations between the outputs of Algorithm \ref{alg LB} are stated in the following theorem.

\vs

\begin{theorem}
Let $\mathcal{A}\in \mathbb{R}^{I_1\times \ldots \times I_N\times J_1\times \ldots \times J_M}$. After $m$ steps of Algorithm \ref{alg LB}, the following relations are satisfied :
\begin{eqnarray}
	\mathcal{A}*_M \mathbb{P}_{m}&=&\mathbb{Q}_m\times_{N+1} B_m^T \label{LB 1}\\
	\mathcal{A}^T *_N \mathbb{Q}_m&=&\mathbb{P}_{m}\times_{M+1}B_m + \mathcal{R}_m \circ e_m, \label{LB 2}
\end{eqnarray}
where $\mathcal{A}*_M\mathbb{P}_m=\left[\mathcal{A}*_M\mathcal{P}_1, \ldots, \mathcal{A}*_M\mathcal{P}_m \right]$ and $\mathcal{A}^T*_N\mathbb{Q}_m=\left[\mathcal{A}^T*_N\mathcal{Q}_1, \ldots, \mathcal{A}^T*_N\mathcal{Q}_m \right]$, with  $\mathcal{P}_1,\ldots, \mathcal{P}$ and $\mathcal{Q}_1,\ldots, \mathcal{Q}_m$ are mutually orthonormal, respectively. The symbol $\circ$ refers to the well-known outer product defined in the literature and $e_m$ is the $m$-th canonical vector of $\mathbb{R}^m$.
\end{theorem}

\vs

\begin{proof}
Let us prove relation \eqref{LB 1}. In order to do so, we first establish that 
\[
\mathbb{Q}_m\times_{N+1} B_m^T = \left[\alpha_1 \mathcal{Q}_1, \beta_1 \mathcal{Q}_1 + \alpha_1 \mathcal{Q}_2 , \ldots, \beta_{m-1}\mathcal{Q}_{m-1}+ \alpha_{m}\mathcal{Q}_m\right] :
\]
By  using the properties of the n-mode product \eqref{n-mode}, we get 
\begin{eqnarray*}
	\left( \mathbb{Q}_m\times_{N+1} B_m^T \right)_{(N+1)}= B_{m}^T \left(\mathcal{Q}_m\right)_{(N+1)}&=& B_{m}^T \begin{bmatrix}
		vec(\mathcal{Q}_1)^T\\
		\vdots \\
		vec(\mathcal{Q}_m)^T\\
	\end{bmatrix}\\
	&=& \begin{bmatrix}
		\alpha_1 vec(\mathcal{Q}_1)^T\\
		\beta_1 vec(\mathcal{Q}_1)^T + \alpha_1 vec(\mathcal{Q}_2)^T \\
		\vdots\\
		\beta_{m-1}vec(\mathcal{Q}_{m-1})^T+ \alpha_{m}vec(\mathcal{Q}_m)^T
	\end{bmatrix},
\end{eqnarray*}
which yields
\[
\mathbb{Q}_m\times_{N+1} B_m^T = \left[\alpha_1 \mathcal{Q}_1, \beta_1 \mathcal{Q}_1 + \alpha_1 \mathcal{Q}_2 , \ldots, \beta_{m-1}\mathcal{Q}_{m-1}+ \alpha_{m}\mathcal{Q}_m\right].
\]
For $j\in \{2,\ldots, m\}$, by from steps $11$ and $16$ of Algorithm \ref{alg LB} , we get 
\[
\alpha_{j}\mathcal{Q}_{j}=\mathcal{A}*_M \mathcal{P}_{j} -\beta_{j-1}\mathcal{Q}_{j-1},
\] 
and, step $1$ gives $\mathcal{A}*_M \mathcal{P}_1=\alpha_{1}\mathcal{Q}_1$. Therefore, we get the relation \eqref{LB 1}.\\

\noindent The proof of Relation \eqref{LB 2} is similar, noticing that 
\[
\mathbb{P}_m\times_{M+1} B_m=\left[ \alpha_1 \mathcal{P}_1 + \beta_1 \mathcal{P}_2, \ldots, \alpha_m \mathcal{P}_m +\beta_m \mathcal{P}_{m+1}, \alpha_m \mathcal{P}_m \right],
\]
and observing that 
\[
\mathcal{R}_m \circ e_m=\left[ 0, \ldots, 0, \mathcal{R}_m\right]\in \mathbb{R}^{J_1\times \ldots \times J_M\times m},
\]
and steps $3$ and $9$ of Algorithm \ref{alg LB}.\\
The orthonormality of $\mathcal{P}_1,\ldots, \mathcal{P}_m$ and $\mathcal{Q}_1,\ldots, \mathcal{Q}_m$ can be proved similarly as for the matrix case.
\end{proof}

In the case when $\mathcal{R}_m$ is normalized, then we can write 
\begin{equation}\label{beta_m}
\mathcal{R}_m=\beta_m \mathcal{P}_{m+1}.
\end{equation}
In that case \eqref{LB 2} can be written as
\[
\mathcal{A}\times_N \mathbb{Q}_m=\mathcal{P}_m\times_{M+1} B_m +\beta_m \mathcal{P}_{m+1}\circ e_m,
\]
or under the following form
\[
\mathcal{A}\times_N \mathbb{Q}_m=\mathbb{P}_{m+1}\times_{M+1} \widetilde{B}_m,
\]
with $\mathbb{P}_{m+1}=\left[\mathbb{P}_m, \mathcal{P}_{m+1}\right]\in \mathbb{R}^{J_1\times \ldots \times J_M \times (m+1)}$ and 
\[
\widetilde{B}_m=\begin{bmatrix}
B_m & \beta_m e_m
\end{bmatrix}\in \mathbb{R}^{m\times (m+1)}.
\]
\subsection{The approximation of first singular triplets}
Let $\mathcal{A}\in \mathbb{R}^{I_1\times \ldots \times I_N \times J_1 \times \ldots \times J_M}$ be a large tensor. Assume that relations \eqref{LB 1} and \eqref{LB 2} are satisfied for  $m\in \mathbb{N}^*$. Consider the svd of  matrix $\mathcal{B}_m$, \textit{i.e.}, 
\[
\mathcal{B}_m=U_m S_m V_m^T=\sum_{i=1}^{m}s_{m,i} u_{m,i} v_{m,i}^T,
\]
where $s_{m,i}=S_m(i,i)$ is the $i$-th singular value of $B_m$, and $u_{m,i}$ and $v_{m,i}$ are the left and the right singular vectors of $B_m$, respectively. $\{s_{m,i},u_{m,i},v_{m,i}\}$ denotes the $i$-th singular triplets of $B_m$. The above decomposition can be written also as
\[
\mathcal{B}_m v_{m,i}= u_{m,i}s_{m,i}, \quad \mathcal{B}_m^T u_{m,i}= v_{m,i}s_{m,i}.
\]
Let $0<k\leq m$, and $\{s_{m,i},u_{m,i},v_{m,i}\}_{i=1}^k$ be the first $k$ singular triplets of $B_m$. The first $k$ approximated singular triplets of $\mathcal{A}$ under the Einstein product $\{\widetilde{s}_{m,i}, \widetilde{\mathcal{U}}_{m,i},\widetilde{\mathcal{V}}_{m,i}\}_{i=1}^k$ are given by 
\begin{equation} \label{approximation}
\widetilde{\mathcal{V}}_{m,i}=\mathbb{P}_m \times_{M+1} v_i^T, \quad \widetilde{\mathcal{U}}_{m,i}=\mathbb{Q}_m \times_{N+1} u_k^T,\quad \widetilde{S}_{m,i}=s_{m,i}, \quad \text{for} \quad i=1,\ldots,k.
\end{equation}

By using the above relations, we notice that for $i=1,\ldots, m$
\begin{eqnarray*}
\mathcal{A}*_M \widetilde{\mathcal{V}}_{m,i}=\mathcal{A}*_M \left(\mathbb{P}_m\times_{M+1} v_{m,i}^T\right)&=&\mathcal{A}*_M\mathbb{P}_m\times_{N+1} v_{m,i}^T\\
&=&\mathbb{Q}_m\times_{N+1}B_m^T\times_{N+1} v_{m,i}^T\\
&=&\mathbb{Q}_m\times_{N+1}\left(v_{m,i}^T B_m^T\right)\\
&=&\mathbb{Q}_m\times_{N+1}\left( B_m v_{m,i}\right)^T\\
&=&\mathbb{Q}_m\times_{N+1}\left( u_{m,i} s_{m,i}\right)^T\\
&=&\mathbb{Q}_m\times_{N+1}  u_{m,i}^T s_{m,i}\\
&=&\widetilde{\mathcal{U}}_{m,i} \widetilde{s}_{m,i}.
\end{eqnarray*}
On the other hand, we have
\begin{eqnarray}
\mathcal{A}^T *_N \widetilde{\mathcal{U}}_{m,i}=\mathcal{A}^T *_N\nonumber \left(\mathcal{Q}_m\times_{N+1}u_{m,i}^T\right)&=&\mathcal{A}^T *_N \mathcal{Q}_m\times_{M+1}u_{m,i}^T \nonumber\\
&=&\left(\mathbb{P}_m\times_{M+1} B_m + \mathcal{R}_m\circ e_m\right)\times_{M+1}u_{m,i}^T\nonumber\\
&=&\mathbb{P}_m\times_{M+1} B_m \times_{M+1}u_{m,i}^T + \left(\mathcal{R}_m\circ e_m\right)\times_{M+1}u_{m,i}^T \nonumber\\
&=&\mathbb{P}_m\times_{M+1} \left( u_{m,i}^T B_m  \right) + \left(\mathcal{R}_m\circ e_m\right)\times_{M+1}u_{m,i}^T \nonumber\\
&=&\mathbb{P}_m\times_{M+1} v_{m,i}^T s_{m,i} + \left(\mathcal{R}_m\circ e_m\right)\times_{M+1}u_{m,i}^T \nonumber\\
&=&\widetilde{\mathcal{V}}_{m,k} \widetilde{s}_{m,i} + \left(\mathcal{R}_m\circ e_m\right) \times_{M+1}u_{m,i}^T \label{eq 3.4}
\end{eqnarray}
Thus, to get better approximation of the first $k^{th}$ singular triplets, it should be found that $\left(\mathcal{R}_m\circ e_m\right) \times_{M+1}u_{m,i}^T$ is small. In other words
\begin{equation}\label{err}
\left\Vert \left(\mathcal{R}_m\circ e_m\right) \times_{M+1}u_{m,i}^T\right\Vert_F=\beta_m \left\vert u_{m,i}(m)\right\vert \leq \varepsilon,
\end{equation}
for some chosen $\varepsilon>0$, where $u_{m,i}(m)$ is the last value of $u_{m,i}$. We have
\begin{eqnarray*}
\left\Vert \left(\mathcal{R}_m\circ e_m\right) \times_{M+1}u_k^T\right\Vert_F &=&\left\Vert \mathcal{R}_m u_{m,i}(m) \right\Vert_F\\
&=&\left\Vert \beta_m \mathcal{P}_m u_{m,i}(m) \right\Vert_F\\
&=& \beta_m \left\vert u_{m,i}(m) \right\vert.
\end{eqnarray*}

\noindent The algorithm summarizes the steps of the algorithm approximating the first $k$ singular triplets of a tensor under the Einstein product.
\begin{algorithm}[H]
\caption{The approximation of the tensor singular triplets by The Einstein Lanczos bidiagonalization  AELB}
\textbf{Input:} $\mathcal{A}\in \mathbb{R}^{I_1\times \ldots \times I_N\times J_1\times \ldots J_M}$, $\mathcal{P}_1\times \in \mathbb{R}^{J_1\times \ldots \times J_N}$ of unit norm, the tensor Lanczos step $m\in \mathbb{N}^*$, the number of the desired singular triplets $k\in \mathbb{N}^*$ such that $k\leq m$.\\
\textbf{Output:} The $k$ first singular triplets $\left\{\widetilde{s}_{m,i}, \widetilde{\mathcal{U}}_{m,i}, \widetilde{\mathcal{V}}_{m,i}\right\}_{i=1}^k$ of $\mathcal{A}$.

\begin{algorithmic}[1]
	\STATE $[\mathbb{P}_m,\mathbb{Q}_m,B_m,\mathcal{R}_m]=ELB(\mathcal{A}, m, \mathcal{P}_1)$.
	\STATE Compute the approximated $i^{th}$ singular triplet $\{\widetilde{s}_{m,i}, \widetilde{\mathcal{U}}_{m,i},\widetilde{\mathcal{V}}_{m,i} \}$ using Eq \ref{approximation}.
\end{algorithmic}
\end{algorithm}

\noindent Notice that as in the case of matrices, the calculation of the last singular triplets (\textit{ie} associated to the smallest singular values) can be made by applying the same method to the inverse of the tensor.
\subsection{Augmentation by Ritz vectors}
Let $\mathcal{A}\in \mathbb{R}^{I_1\times \ldots\times I_N  \times J_1 \times \ldots \times J_M}$ be a large tensor. Assuming that relations \eqref{LB 1} and \eqref{LB 2} are satisfied for a certain $m\leq \min(I_1\ldots I_N, J_1\ldots J_M)$, we would like to approximate the first $k$ singular triplets of $\mathcal{A}$ with $k<m$.\\

\noindent The approximated right singular tensor $\widetilde{\mathcal{V}}_{m,k}$ is a Ritz tensor of $\mathcal{A}^T*_N\mathcal{A}$ associated with the Ritz value $\widetilde{s}_{i,m}^2$, and verifies:
\[
\mathcal{A}^T*_N\mathcal{A}*_M\widetilde{\mathcal{V}}_{m,i}=\widetilde{\mathcal{V}}_{i,m} \widetilde{s}_{m,i} + \left(\mathcal{R}_m \circ e_m\right)\times_{M+1}u_i^T\widetilde{s}_{m,i}.
\]

\noindent Assume that $\beta_m$ in \eqref{beta_m} is nonvanishing (because in the opposite cas, the singular triplets are well approximated, for more details see \cite{elha4,baglama2005augmented}).\\
The key idea of augmentation by Ritz tensors is reconstructing similar equations as \eqref{LB 1} and \eqref{LB 2}, only changing the projection subspaces. Instead of choosing the first $k$ tensors of $\mathbb{P}_m$, namely $\mathcal{P}_1,\ldots, \mathcal{P}_k$, as shown in Algorithm \ref{alg LB}, we replace them by the first $k$ Ritz tensors.\\

\vs

\noindent Let us consider the following tensor :

\vs

\begin{equation}\label{Ritz 1}
\widecheck{\mathbb{P}}_{k+1}=\left[\widetilde{\mathcal{V}}_{m,1},\ldots, \widetilde{\mathcal{V}}_{m,k},\mathcal{P}_{m+1}\right]\in \mathbb{R}^{J_1\times \ldots \times J_M\times (k+1)}.
\end{equation}
It can be seen that 
\[
\mathcal{A}*_M \widecheck{\mathbb{P}}_{k+1}=\left[\widetilde{\mathcal{U}}_{m,1}\widetilde{s}_{1,m},\ldots, \widetilde{\mathcal{U}}_{m,k}\widetilde{s}_{m,k}, \mathcal{A}*_N \mathcal{P}_{m+1}\right].
\]
Orthogonalizing $\mathcal{A}*_M \mathcal{P}_{m+1}$ against $\left\{\widetilde{\mathcal{U}}_{m,i}\right\}_{i=1}^k$ gives
\[
\mathcal{A}*_M \mathcal{P}_{m+1}=\sum_{i=1}^{k}\rho_i \widetilde{\mathcal{U}}_{m,i}+\widecheck{\mathcal{R}}_m,
\]
with 
\begin{eqnarray*}
\rho_i=\langle \widetilde{\mathcal{U}}_{m,i}, \mathcal{A}*_N \mathcal{P}_{m+1} \rangle &=& {\tt tr}\left(\widetilde{\mathcal{U}}_{i,m}^T*_N \mathcal{A}*_N \mathcal{P}_{m+1}\right)\\
&=&\beta_m \langle \left(\mathcal{P}_{m+1}\circ e_m\right)\times_{M+1} u_{m,i}^T, \mathcal{P}_{m+1} \rangle\\
&=&\beta_m \langle u_{m,i}, e_m \rangle=\beta_m u_{m,i}(m).
\end{eqnarray*}
Normalization of $\widecheck{\mathcal{R}}_k$, gives $\widecheck{\mathcal{R}}_k=\widecheck{\alpha}_{k+1}\widecheck{\mathcal{R}}_k'$. Let us now consider the following tensor and matrix
\begin{equation}\label{Ritz 2}
\widecheck{\mathcal{Q}}_{k+1}=\left[  \widetilde{\mathcal{U}}_{m,1},\ldots, \widetilde{\mathcal{U}}_{m,k}, \widecheck{\mathcal{R}}_k'  \right]\in \mathbb{R}^{I_1 \times \ldots \times I_N \times (k+1)},
\end{equation}
and 
\begin{equation}\label{Ritz 3}
\widecheck{B}_{k+1}=\begin{pmatrix}
	\widetilde{s}_{m,1} & & & \rho_1\\
	& \ddots & & \vdots\\
	& & \widetilde{s}_{m,k}& \rho_k\\
	& & & \widecheck{\alpha}_{k+1}. 
\end{pmatrix}\in \mathbb{R}^{(k+1)\times (k+1)}.
\end{equation}
Since we have 
\[
\mathcal{A}*_M\widecheck{\mathbb{P}}_{k+1}=\left[ \widetilde{\mathcal{U}}_{m,1}\widetilde{s}_{m,1}, \ldots, \widetilde{\mathcal{U}}_{m,k} \widetilde{s}_{k,m}, \sum_{i=1}^{k} \rho_i \widetilde{\mathcal{U}}_{m,i} + \widecheck{\mathcal{R}}_k  \right],
\]
\noindent we obtain
\begin{equation}\label{eq Ritz 1}
\mathcal{A}*_M \widecheck{\mathbb{P}}_{k+1}=\widecheck{\mathbb{Q}}_{k+1}\times_{N+1}\widecheck{B}_{k+1}^T.
\end{equation}

\noindent On the other hand, we have
\[
\mathcal{A}^T *_N \widecheck{\mathcal{Q}}_{k+1}=\left[  \mathcal{A}^T*_N \widetilde{\mathcal{U}}_{m,1}, \ldots,   \mathcal{A}^T*_N \widetilde{\mathcal{U}}_{m,k},  \mathcal{A}^T*_N \widecheck{\mathcal{R}}_{k}' \right].
\]

\noindent Moreover, from \eqref{eq 3.4} we get
\[
\mathcal{A}^T*_N \widetilde{\mathcal{U}}_{m,i}=\widetilde{\mathcal{V}}_{m,i} \widetilde{s}_{m,i} + \mathcal{P}_{m+1}\rho_i.
\]

\noindent It can be easily shown  that $\left(\mathcal{A}^T *_N \widecheck{\mathcal{R}}_k'\right)$ is orthogonal to $\widetilde{\mathcal{V}}_{m,i}$ (by showing that $\langle \mathcal{A}^T *_N \widecheck{\mathcal{R}}_k', \widetilde{\mathcal{V}}_{m,i} \rangle=0$). Thus we can write :
\[
\mathcal{A}^T *_N \widecheck{\mathcal{R}}_k'=\gamma \mathcal{P}_{m+1} + \mathcal{F}_{k+1},
\]
\noindent where $\mathcal{F}_{k+1}$ is orthogonal to $\mathcal{P}_{m+1}$ as well as $\widetilde{\mathcal{V}}_{m,i}$, with $\gamma=\widecheck{\alpha}_{k+1}$.\\
Consequently, we have
\begin{eqnarray*}
\mathcal{A}^T*_N \widecheck{\mathcal{Q}}_{k+1}&=&\left[  \widetilde{\mathcal{V}}_{m,1}\widetilde{s}_{m,1} + \mathcal{P}_{m+1} \rho_1, \ldots, \widetilde{\mathcal{V}}_{m,k}\widetilde{s}_{m,k} + \mathcal{P}_{m+1} \rho_k, \widecheck{\alpha}_{k+1}\mathcal{P}_{m+1} + \mathcal{F}_{k+1} \right]\\
&=& \widecheck{\mathbb{P}}_{k+1}\times_{M+1}\widecheck{B}_{k+1} + \mathcal{F}_{k+1}\circ e_{k+1}.
\end{eqnarray*}
By normalizing $\mathcal{F}_{k+1}$, we get $\mathcal{F}_{k+1}=\widecheck{\beta}_{k+1}\widecheck{\mathcal{P}}_{k+2}$. Therefore
\begin{equation}\label{eq Ritz 2}
\mathcal{A}^T*_N \widecheck{\mathcal{Q}}_{k+1}=\widecheck{\mathbb{P}}_{k+1}\times_{M+1}\widecheck{B}_{k+1} + \widecheck{\beta}_{k+1}\widecheck{\mathcal{P}}_{k+2}\circ e_{k+1}.
\end{equation}

\noindent If $\widecheck{\beta}_{k+1}=0$, then the singular triplets are well approximated. Let us assume the opposite. In order to obtain similar relations as \eqref{LB 1} and \eqref{LB 2}, we may append new tensors to $\widecheck{\mathbb{P}}_{k+1}$ and $\widecheck{\mathbb{Q}}_{k+1}$. Let us assume that
\begin{eqnarray*}
\widecheck{\mathbb{P}}_{k+1}&=&\left[\widecheck{\mathcal{P}}_{1},\widecheck{\mathcal{P}}_2, \ldots, \widecheck{\mathcal{P}}_{k+1}\right],\\
\widecheck{\mathbb{Q}}_{k+1}&=&\left[\widecheck{\mathcal{Q}}_{1},\widecheck{\mathcal{Q}}_2, \ldots, \widecheck{\mathcal{Q}}_{k+1}\right],
\end{eqnarray*}
and 
\[
\widecheck{\mathbb{P}}_{k+2}=\left[\widecheck{\mathbb{P}}_{k+1}, \widecheck{\mathcal{P}}_{k+2}\right]\in \mathbb{R}^{J_1\times \ldots \times J_M \times (k+2)}.
\]
Let $\left(\mathcal{I}_N - \widecheck{\mathcal{Q}}_{k+1}*_N \widecheck{\mathcal{Q}}_{k+1}^T\right)*_N \mathcal{A}*_M \widecheck{\mathcal{P}}_{k+2}$ be the orthogonalization of $\mathcal{A}*_M \widecheck{\mathcal{P}}_{k+2}$ against $\widecheck{\mathcal{Q}}_1,\ldots, \widecheck{\mathcal{Q}}_{k+1}$. Its normalization gives
\[
\widecheck{\alpha}_{k+2}\widecheck{\mathcal{Q}}_{k+2}=\left(\mathcal{I}_N - \widecheck{\mathcal{Q}}_{k+1}*_N \widecheck{\mathcal{Q}}_{k+1}^T\right)*_N \mathcal{A}*_M \widecheck{\mathcal{P}}_{k+2}.
\] 
Simple calculations yield
\[
\widecheck{\alpha}_{k+2}\widecheck{\mathcal{Q}}_{k+2}=\mathcal{A}*_M \widecheck{\mathcal{P}}_{k+2}- \widecheck{\beta}_{k+1}\widecheck{\mathcal{Q}}_{k+1}.
\]
\noindent  Equation \eqref{eq Ritz 1}, we get 
\begin{equation}\label{eq Ritz 3}
\mathcal{A}*_M \widecheck{\mathbb{P}}_{k+2}=\widecheck{\mathbb{Q}}_{k+2}\times_N \widecheck{B}_{k+2}^T,
\end{equation}
with $\widecheck{\mathbb{Q}}_{k+2}=\left[\widecheck{\mathbb{Q}}_{k+1}, \widecheck{\mathcal{Q}}_{k+2}\right]\in \mathbb{R}^{I_1\times \ldots \times I_N \times (k+2)}$ and 
\[
\widecheck{B}_{k+2}=\begin{pmatrix}
\widetilde{s}_{m,1} & & & \rho_1 & 0\\
&\ddots & & \vdots & \vdots\\
& & & &\\
& & \widetilde{s}_{m,k} & \rho_k & 0\\
& & & \widecheck{\alpha}_{k+1} & \widecheck{\beta}_{k+1}\\
& & & & \widecheck{\alpha}_{k+2}
\end{pmatrix}\in \mathbb{R}^{(k+2)\times (k+2)}.
\]
On the other hand, a new tensor $\widecheck{\mathcal{P}}_{k+3}$ of $\widecheck{\mathbb{P}}_{k+3}$ can be obtained by orthonormalizing $\mathcal{A}^T *_N \widecheck{\mathcal{Q}}_{k+2}$ against $\widecheck{\mathcal{P}}_{1}, \ldots, \widecheck{\mathcal{P}}_{k+2}$, \textit{i.e.},
\[
\widecheck{\beta}_{k+2}\widecheck{\mathcal{P}}_{k+3}=\left(\mathcal{I}_M - \widecheck{\mathcal{P}}_{k+2}*_N \widecheck{\mathcal{P}}_{k+2}^T \right)*_M \mathcal{A}^T *_N \widecheck{\mathcal{Q}}_{k+2},
\] 
giving
\[
\widecheck{\beta}_{k+2}\widecheck{\mathcal{P}}_{k+3}=\mathcal{A}^T *_N \widecheck{\mathcal{Q}}_{k+2} - \widecheck{\beta}_{k+2}\widecheck{\mathcal{P}}_{k+3}\circ e_{k+2}.
\] 
Therefore, we have
\begin{equation}
\mathcal{A}^T *_N \widecheck{\mathcal{Q}}_{k+2}=\widecheck{\mathcal{P}}_{k+2}\times_{M+1} \widecheck{B}_{k+2} + \widecheck{\beta}_{k+2}\widecheck{\mathcal{P}}_{k+3}\circ e_{k+2}.\label{eq Ritz 4}
\end{equation}
\noindent Repeating the above process $m-k+2$ times in order to compute $\widecheck{\mathcal{Q}}_{k+2}$ and $\widecheck{\mathcal{P}}_{k+3}$, we find the following decompositions which are analogous to \eqref{LB 1} and \eqref{LB 2} :
\begin{eqnarray}
\mathcal{A}*_M \widecheck{\mathbb{P}}_m&=&\widecheck{\mathbb{Q}}_{m}\times_{N+1} \widecheck{B}_{m}^T, \label{Ritz 4}\\
\mathcal{A}^T *_N \widecheck{\mathbb{Q}}_m&=&\widecheck{\mathbb{P}}_{m}\times_{M+1} \widecheck{B}_{m}+ \widecheck{\beta}_{m} \widecheck{\mathcal{P}}_{m+1}\circ e_m, \label{Ritz 5}
\end{eqnarray}
with $\widecheck{\mathbb{P}}_m \in \mathbb{R}^{J_1\times \ldots \times J_N \times m}$ and $\widecheck{\mathbb{Q}}_m \in \mathbb{R}^{I_1\times \ldots \times I_N \times m}$ having their sub-tensors orthonormal, \textit{i.e.}, $\widecheck{\mathcal{P}}_1,\ldots, \widecheck{\mathcal{P}}_m$ and $\widecheck{\mathcal{Q}}_1, \ldots, \widecheck{\mathcal{Q}}_m$ are mutually orthonormal, respectively. Moreover, $\widecheck{\mathcal{P}}_{m+1}$ is orthogonal to all the sub-tensors of $\widecheck{\mathcal{P}}_{m}$; $\widecheck{\mathcal{P}}_1,\ldots, \widecheck{\mathcal{P}}_m$, and $\widecheck{B}_{m}$ is an upper triangular matrix given by 
\begin{equation}\label{Ritz 6}
\widecheck{\mathcal{B}}_{m}=\begin{pmatrix}
	\widetilde{s}_{m,i} & & & \rho_1 & & & 0\\
	& \ddots & &  \vdots & & & \vdots \\
	& & \widetilde{s}_{m,k} & \rho_k & & & \\
	& & & \widecheck{\alpha}_{k+1} & \widecheck{\beta}_{k+1} & & \\
	& & & & \ddots & \ddots & \\
	& & & & & \widecheck{\alpha}_m & \widecheck{\beta}_{m-1}\\
	& & & & & & \widecheck{\alpha}_{m}
\end{pmatrix}\in \mathbb{R}^{m\times m}.
\end{equation}
The new residual is then 
\begin{equation}\label{neq res}
\widecheck{\mathcal{R}}_m=\widecheck{\beta}_m \widecheck{\mathcal{P}}_{m+1}.
\end{equation}

\noindent The following algorithm summarises the process of Ritz augmentation

\begin{algorithm}[H]
\caption{Ritz augmentation for the tensor approximated SVD based on the Einstein product (LBR)}
\label{alg Ritz}
\textbf{Input:}
\begin{itemize}
	\item $\mathcal{A}\in \mathbb{R}^{I_1\times\ldots \times I_N\times J_1 \times \ldots \times J_M}$
	\item $\mathcal{P}_1\in \mathbb{R}^{J_1\times \ldots \times J_M}$ of unit norm
	\item the Lanczos bidiagonalization steps $m$
	\item  the number of the desired singular triplets $k$
	\item  $\varepsilon >0$
\end{itemize}

\textbf{Output:} The largest $k$ singular triplets $\left\{ \widetilde{s}_{m,i}, \widetilde{\mathcal{U}}_{m,i},\widetilde{\mathcal{V}}_{m,i}  \right\}_{i=1}^k$.

\begin{algorithmic}[1]
	\STATE $\left[ \mathbb{P}_m, \mathbb{Q}_m, B_m, \mathcal{R}_m\right]=RLB(\mathcal{A}, m,\mathcal{P}_1)$.
	\STATE Compute the SVD of $B_m$; $[U_m,S_m,V_m]={\tt svd}\left(B_m\right)$.
	\STATE Check the convergence of the first $k$ singular triplets using Equation \eqref{approximation}
	\IF {They all converge}
	\STATE Exit
	\ELSE
	\STATE Compute $\widecheck{\mathbb{P}}_{k+1}$, $\widecheck{\mathbb{Q}}_{k+1}$, and  $\widecheck{B}_{k+1}$, from \eqref{Ritz 1}, \eqref{Ritz 2}, \eqref{Ritz 3}, respectively.
	\STATE Append new tensors to $\widecheck{\mathbb{P}}_{k+1}$ and $\widecheck{\mathbb{Q}}_{k+1}$ to construct  $\widecheck{\mathbb{P}}_{m}$, $\widecheck{\mathbb{Q}}_{m}$, and  $\widecheck{B}_{m}$ given in \eqref{Ritz 4} ,\eqref{Ritz 5}, and \eqref{Ritz 6}, and a new residual $\widecheck{\mathcal{R}}_m$ given in \eqref{neq res}. 
	\STATE Consider $\mathbb{P}_m=\widecheck{\mathbb{P}}_m$, $\mathbb{Q}_m=\widecheck{\mathbb{Q}}_m$, $B_m=\widecheck{B}_m$, $\mathcal{R}_m=\widecheck{\mathcal{R}}_m$.
	\STATE Go to $2$.
	\ENDIF
\end{algorithmic}
\end{algorithm}

\noindent In the case we would like to approximate the smallest singular triplets of $\mathcal{A}$, the same process can be applied using the last $k$ Ritz tensors instead of the first ones.
\section{Principal Component Analysis based on the Einstein product}
\label{sec PCA}

\noindent Principal Component Analysis (PCA) is a very well-established technique in statistics and engineering applications such as in data denoising, classification problems, and facial recognition \cite{Hached,elha4,class1,class2}. In the case of image classification and facial recognition, color images are represented as tensors which leads, in most cases, to a large  training tensor. The adaption of PCA process and Krylov subspaces methods to image recognition was first handled in \cite{Hached}, using the t- or c-product in the case of 3rd-order tensors (RGB images). Then, in \cite{elha4}, the authors proposed computational refinements to this method for better efficiency. This section presents a PCA process based on the Einstein product.

\vs

\noindent Assume we have $N$ training color images of the same size. Each image is represented by a third-order tensor $\mathcal{I}_1,\ldots, \mathcal{I}_N$ of size $\ell \times p \times n$.

\vs 

\noindent The steps of the PCA algorithm for  facial recognition are listed as follows:

\vs

\begin{itemize}
\item Vectorize each of the training images, such that each vector is going to be of size $\ell p\times n$. This vectorization is implemented in Matlab by the command Matlab $X_i={\tt reshape}\left(\mathcal{I}_i,[\ell p, n]\right)$.
\item Compute the mean image vector of $X_1,\ldots, X_N$, \textit{i.e.},
\[
M=\dfrac{\displaystyle{\sum_{i=1}^{N}}X_i}{N}.
\]
\item Construct the  tensor $\widebar{\mathcal{X}}$, such that $\widebar{\mathcal{X}}(:,:,i)=X_i-M$.
\item Compute the first $k$ left tensor singular vectors $U_1,\ldots, U_k \in \mathbb{R}^{\ell \times p }$. Construct the projector tensor 
\[
\mathcal{U}_k=\left[  U_1,\ldots, U_k \right]\in \mathbb{R}^{\ell \times p\times k}.
\]
\item Project all the training images onto the the space $span\{U_1, \ldots, U_k\}$ (the projection space)  by using the formula $\mathcal{U}_k^T*_2 \widebar{\mathcal{X}}.$
\item Project a test image represented by the tensor $\mathcal{I}_0\in \mathbb{R}^{\ell \times p \times n}$, by the formula $\mathcal{U}_k^T*\left( {\tt reshape}\left(\mathcal{I}_0,[\ell p, n]\right)- M  \right)$.
\item Computed the smallest distance between the projected training images and the test image.
\end{itemize}

The Algorithm describes the PCA under the Einstein product for the facial recognition as follows.
\begin{algorithm}[H]
\caption{Facial recognition based of the Einstein product's PCA}\label{alg 3}
\textbf{Inputs:}
\begin{itemize}
	\item $I_1,\ldots, I_N$ training images
	\item $I_0$ test image
	\item number $k>0$ : the desired number of singular triplets
	
\end{itemize}
\textbf{Outputs:} Closest image in the training database.
\begin{algorithmic}[1]
	\STATE Construct $M$ and $\widebar{\mathcal{X}}$.
	\STATE $[\mathcal{U}_k,\mathcal{S},\mathcal{V}]=AELB(\widebar{X},\mathcal{P}_1, m,k)$  or  $[\mathcal{U}_k,\mathcal{S},\mathcal{V}]=LBR(\widebar{X},\mathcal{P}_1, m,k,\varepsilon)$.
	\STATE Compute the projected training images $\mathcal{P}=\mathcal{U}_k^T* \left( {\tt reshape}\left(\mathcal{I}_0,[\ell p, n]\right)- M  \right)$.
	\STATE Find $d=\min_{i=1,\ldots, N}\left\Vert   \mathcal{P}(:,:,i)-\mathcal{P}_0  \right\Vert_F.$
\end{algorithmic}
\end{algorithm}

\noindent In line $2$ of the above algorithm, we have the choice to use either Lanczos bidiagonalization process or Ritz augmentation to approximate the first $k$ singular triplets.

\section{Numerical experiments}
\label{sec 5}
In this section, we present some numerical results and applications related to the proposed methods. This section is divided into three subsections: the first one is devoted to some tests of Algorithms \ref{alg LB} and \ref{alg Ritz} on synthetic data. In the second one, Algorithms \ref{alg LB} and \ref{alg Ritz} are applied to data compression. The results of the application of the algorithms on facial recognition are depicted in the last subsection.
All computations are carried out on a laptop computer with 2.3 GHz Intel Core i5 processors and 8 GB of memory using MATLAB 2018a.

\subsection{Synthetic data}
In this part, tests of Algorithms \ref{alg LB} (LB) and \ref{alg Ritz} (Ritz) are performed on synthetic data. the first sub-part is devoted  to the LB method to approximate the first $k$ singular triplets, the second one is reserved for the results of the method of Ritz to approximate the first and the last singular triplets. The synthetic data is obtained by the Matlab command ${\tt randn}$.
Two factors are used to measure the quality of the approximation of the first $k$ singular triplets. Let $\{ s_i, \mathcal{U}_i, \mathcal{V}_i\}_{i=1}^k$ be the first $k$ singular triplets of a tensor $\mathcal{A}\in\mathbb{R}^{I_1\times \ldots I_N \times I_1 \times \ldots\times I_N}$. The residual norm defined by 
\[
\text{Res.norm}=\left\Vert \mathcal{A}*_N \mathcal{V}_i - \mathcal{U}_i s_i \right\Vert_F, \quad i=1,\ldots,k.
\]
Also, the Global Residual norm is defined by 
\[
\text{GRes.norm}=\left\Vert \mathcal{A}*_N \mathcal{V} - \mathcal{U}*_N \mathcal{S} \right\Vert_F,
\]
where $\mathcal{U},\mathcal{V}$, and $\mathcal{S}$ are constructed from the triplets $\{ s_i, \mathcal{U}_i, \mathcal{V}_i\}$ for $i=1,\ldots, k$.

\vs

\subsubsection{Lanczos bidiagonalization}

\noindent In Table \ref{tab 1}, the value of the residual norm, denoted \textbf{Res.norm} 
of each approximated singular triplet is shown, for various sizes of tensor $\mathcal{A}$. The number of singular triplets to be approximated is equal to the number of the Tensor Lanczos bidiagonalization steps, chosen to be equal to four in this example. 
\begin{table}[H]
\centering
\begin{tabular}{l | c | c  |c}
	\hline size($\mathcal{A}$) &$50\times 20 \times 50 \times 20$ & $50 \times 100 \times 50 \times 100$ & $100\times 50 \times 20 \times 100\times 50 \times 20$\\
	\hline
	$i=1$ & 3.319e-14    &3.483e-13    & 7.097e-13    \\
	$i=2$ & 4.163e-14    &3.154e-13   & 6.760e-13    \\
	$i=3$ & 2.949e-14     &2.668e-13     & 5.356e-13     \\
	$i=4$ & 2.444e-14 &2.249e-13  & 4.482e-13  \\
	\hline 
\end{tabular}
\captionof{table}{The Res.norm of each approximated singular triplet when $k=4$ and $m=4$ for tensors of different sizes.}\label{tab 1}
\end{table}

\noindent In Table \ref{tab 2} we computed the \textbf{GRes.norm} for the same experiment.

\begin{table}[H]
\centering
\begin{tabular}{l | c | c  |c}
	\hline size($\mathcal{A}$) &$50\times 20 \times 50 \times 20$& $50 \times 100 \times 50 \times 100$& $100\times 50 \times 20 \times 100\times 50 \times 20$\\
	\hline GRes.norm   &1.3145e-13 & 8.5752e-13& 3.5350e-12  \\
	\hline  
\end{tabular}
\captionof{table}{The GRes.norm of each approximated singular triplet when $k=4$ and $m=4$ for tensors of different sizes.}\label{tab 2}
\end{table}

\noindent Tables \ref{tab 1} and \ref{tab 2} show the effectiveness of the presented method in terms of accuracy. In Table \ref{tab 3}, we give the CPU-time of the suggested method against using the Exact Einstein SVD for different sizes of the used tensors. We chose to approximate the $k=5$ largest singular triplets.

\begin{table}[H]
\centering
\begin{tabular}{l|c|c|c}
	\hline size($\mathcal{A}$)& $50\times 20 \times 50 \times 20$& $50 \times 100 \times 50 \times 1000$& $100\times 50 \times 20 \times 100\times 50 \times 20$\\
	\hline
	Exact & 0.7402 &52.9240 & 574.5531 \\
	Approximated & 0.4796& 3.4324& 7.0692   \\
	\hline 
\end{tabular}
\captionof{table}{CPU-time of the approximated method vs exact computation, for $k=5$ and $m=5$.}\label{tab 3}
\end{table}

\begin{figure}[H]
\centering
\begin{tabular}{c}
	\includegraphics[width=0.5\linewidth]{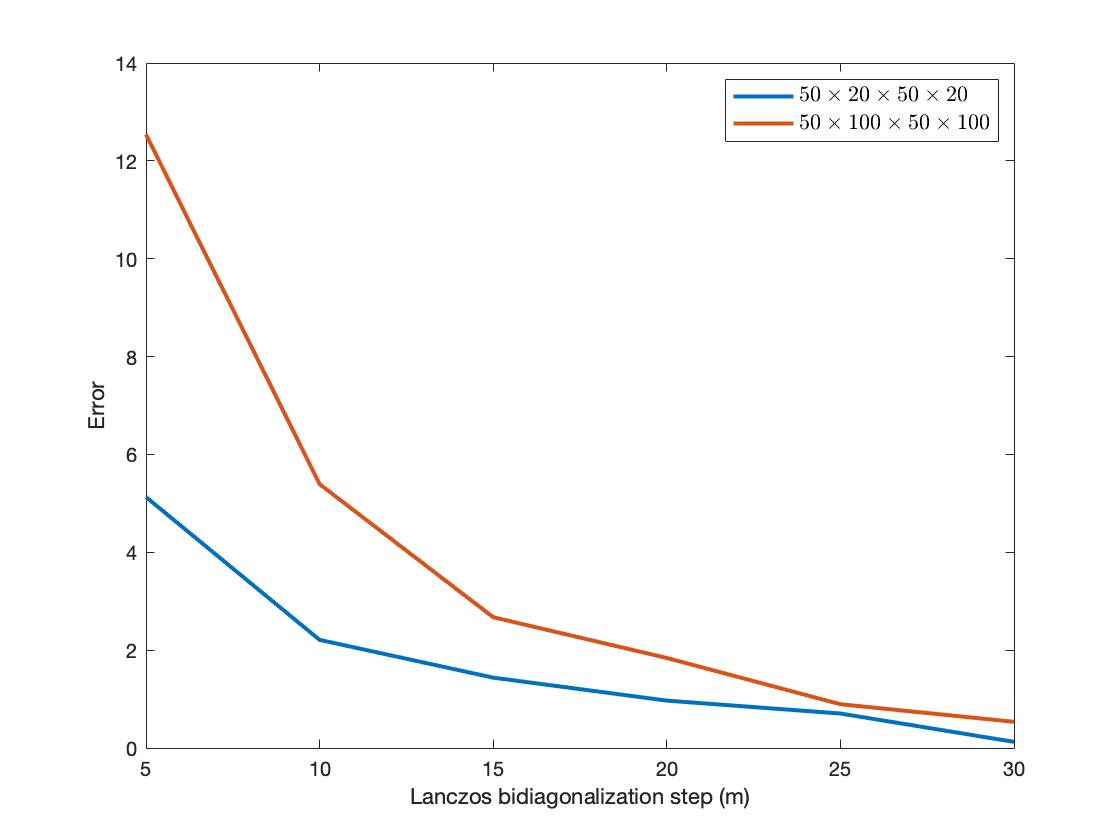}
\end{tabular}
\captionof{figure}{The error norm $\beta_m\left\vert u_1(m)\right\vert$ for different values of $m$.}\label{fig 1}
\end{figure}

\noindent The figures displayed in Table \ref{tab 3} show that the Lanczos bidiagonalization method represents an important gain in terms of CPU-time, which becomes considerable as the size of the tensor gets larger. \\
\noindent In Figure \ref{fig 1}, we plot the error norm as defined in \ref{err} in function of the number $m$  of Lanczos bidiagonalization steps when computing the first singular triplet, for synthetic tensors of sizes $50\times 20\times 50 \times 20$ and $50\times 100\times 50 \times 100$, respectively.

\noindent From Figure \ref{fig 1}, it can be observed that when the number $m$ of Lanczos bidiagonalization gets larger, the error decreases. 

\subsubsection{Ritz augmentation.}
In this part, the numerical results of some tests on synthetic data are illustrated. These tests are carried out by using the  Ritz augmentation process to approximate a desired number of  tensor singular triplets (this number is denoted by $k$), for some chosen  values of Lanczos bidiagonalization step $m>k$. In the sequel, the $i$-th approximated singular value is denoted by $\widetilde{s}_i$, while the exact one is denoted by $s_i$.

\noindent \textbf{Approximation of the first singular triplets}\\
\noindent This part is devoted to the approximation of the first (largest) singular triplets of synthetic tensors using the Ritz augmentation process.\\
In Table \ref{tab ritz1}, the values of Res.norm are given when using the Ritz process to approximate the first four singular triplets of tensors for different tensor sizes, and by using $m=15$. Table \ref{tab ritz2} shows the error between the approximated singular values and the exact ones.
\begin{table}[H]
\centering
\small\addtolength{\tabcolsep}{-3pt}
\begin{tabular}{l|c|c|c|c}
	\hline size($\mathcal{A}$) & $100\times 100\times 50$& $50\times 20\times 50\times 20$ & $50\times 100\times 50\times 100$& $ 50\times 20\times 10\times 50 \times 20\times 10$ \\
	\hline
	$i=1$&1.77e-13 & 1.54e-13&6.55e-13 & 3.56e-12  \\
	$i=2$& 2.98e-13&1.76e-13 &1.95e-12& 2.16e-12  \\
	$i=3$&2.26e-13&1.99e-13 & 1.74e-12&  1.28e-12\\
	$i=4$ &1.49e-13&2.21e-13 &4.90e-13& 1.19e-12 \\
	\hline 
\end{tabular}
\captionof{table}{The Res.norm of each of the approximated singular triplets when using Ritz augmentation for $m=15$ and $k=4$.}
\label{tab ritz1}
\end{table}

\begin{table}[H]
\centering
\small\addtolength{\tabcolsep}{-3pt}
\begin{tabular}{l|c|c|c|c}
	\hline size($\mathcal{A}$)& $100\times 100\times 50$& $50\times 20\times 50\times 20$ & $50\times 100\times 50\times 100$& $ 50\times 20\times 10\times 50 \times 20\times 10$ \\
	\hline
	$i=1$&2.84e-14& 2.13e-13&5.68e-14 &4.26e-13  \\
	$i=2$& 3.12e-13&1.98e-13 &8.52e-13&7.95e-13  \\
	$i=3$&1.98e-13 &9.94e-14 &1.42e-13& 5.40e-13  \\
	$i=4$ &8.52e-14&7.64e-11 &2.84e-10&  1.17e-09 \\
	\hline 
\end{tabular}
\captionof{table}{The values of $\left\vert \widetilde{s}_i-s_i\right\vert$  of each of the approximated singular triplets when using Ritz augmentation for $m=15$ and $k=4$.}
\label{tab ritz2}
\end{table}
\noindent Table \ref{tab ritz3}, presents the GRes.norm for different tensors when the first four singular triplets are approximated using the Ritz process with $m=15$. Tables \ref{tab ritz4} shows the influence of $m$ (Lanczos bidiagonalization step) on the GRes.norm and the number of iterations (Iter), when the first four singular triplets are desired to be approximated.

\begin{table}[H]
\centering
\small\addtolength{\tabcolsep}{-3pt}
\begin{tabular}{l|c|c|c|c}
	\hline
	size($\mathcal{A}$) & $100\times 100\times 50$& $50\times 20\times 50\times 20$ & $50\times 100\times 50\times 100$& $ 50\times 20\times 10\times 50 \times 20\times 10$  \\
	\hline
	GRes.norm &4.40e-13 &7.57e-13 &5.37e-12 & 1.75e-11  \\
	\hline
\end{tabular}
\captionof{table}{The GRes.norm of the approximated first singular triplets when using Ritz augmentation for $m=15$ and $k=4$.}
\label{tab ritz3}
\end{table}

\begin{table}[H]
\small\addtolength{\tabcolsep}{-2pt}
\centering
\begin{tabular}{l|c|c|c|c|c|c}
	\hline
	size($\mathcal{A}$) & \multicolumn{2}{l|}{$100\times 100\times 50$}& \multicolumn{2}{l|}{$50\times 20\times 50\times 20$} & \multicolumn{2}{|l}{$50\times 100\times 50\times 100$}  \\
	\hline & GRes & Iter& GRes & Iter& GRes & Iter\\
	\hline
	$m=5$&1.63e-12 &107 &4.63e-11  &1000 &1.84e-10&1000\\
	$m=10$&9.40e-13  &13 &2.98e-12& 45& 9.99e-12&86 \\
	$m=15$& 4.40e-13 &5 &7.57e-13 &11 &5.37e-12 &19\\
	$m=20$& 7.85e-13  &3 &8.54e-13 &6 &4.08e-12  &11\\
	\hline
\end{tabular}
\captionof{table}{The values of the GRes.norm and the number of iterations, when using Ritz augmentation with $k=4$ and $m=5,10,15,20$.}
\label{tab ritz4}
\end{table}

\begin{figure}[H]
\centering
\begin{tabular}{c}
	\includegraphics[width=0.5\linewidth]{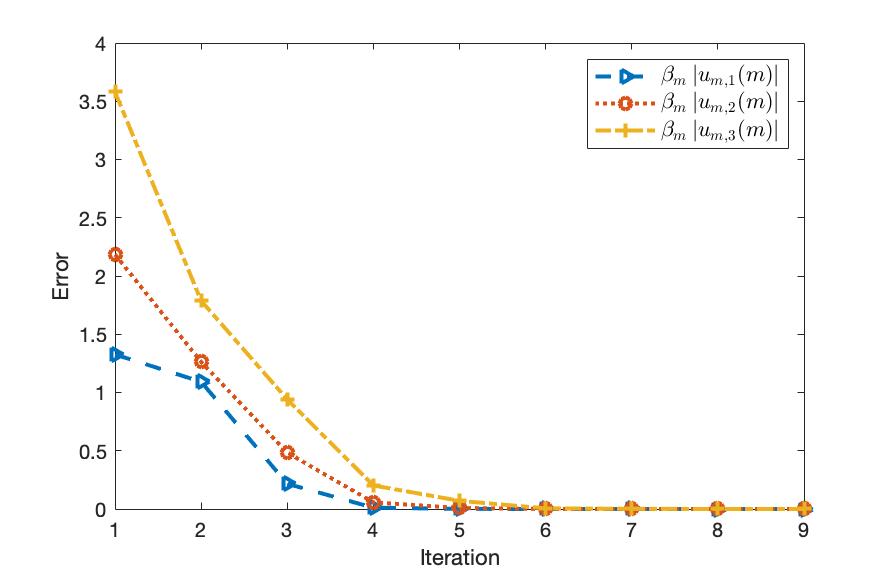}
\end{tabular}
\captionof{figure}{The error norm on each of the approximated first three singular triplets of the tensor of size $50\times20\times 50\times 20$ when $m=15$.}\label{fig 5.2}
\end{figure}

\noindent Figure \ref{fig 5.2} illustrates the error norm as defined in \eqref{err} on the first three singular triplets of the tensor of size $50\times 20\times 50 \times 20$ for $m=15$.

\noindent Figure \ref{fig 5.2} show that the considered singular values are accurately approximated after just a few iterations.

\noindent \textbf{Approximation of the last singular triplets}

\noindent  In this part, let us assume that a tensor admits $k_{max}$ singular triplets. Table \ref{tab last1} shows the Res.norm when approximating the last four singular triplets with $m=15$. Table \ref{tab last2} represents the error between the approximated singular triplets and the exact ones, when $m=15$ and $k_{max}=4$.
\begin{table}[H]
\centering
\small\addtolength{\tabcolsep}{-3pt}
\begin{tabular}{l|c|c|c}
	\hline size($\mathcal{A}$) & $100\times 100\times 50$& $20\times 10\times 20\times 10$ & $50\times 20\times 50\times 20$  \\
	\hline $ k=k_{max}-3$& 2.37e-13 &   5.24e-14 & 2.91e-12 \\
	$ k=k_{max}-2$&2.11e-13& 6.75e-14 & 2.66e-12  \\
	$ k=k_{max}-1$&1.69e-13&  4.64e-14&1.60e-12 \\
	$ k=k_{max}$ &1.02e-13& 2.76e-14&1.26e-12  \\
	\hline 
\end{tabular}
\captionof{table}{The Res.norm of each of the approximated last singular triplets when using Ritz augmentation for $m=15$ and $k_{max}=4$.}
\label{tab last1}
\end{table}

\begin{table}[H]
\centering
\small\addtolength{\tabcolsep}{-3pt}
\begin{tabular}{l|c|c|c}
	\hline
	size($\mathcal{A}$) & $100\times 100\times 50$& $20\times 10\times 20\times 10$& $50\times 20\times 50\times 20$  \\
	\hline $ k=k_{max}-3$& 7.10e-14& 1.12e-12    &1.18e-10    \\
	$ k=k_{max}-2$& 2.13e-13 &  1.11e-13  & 2.11e-13  \\
	$ k=k_{max}-1$& $<$machine eps.& 8.30e-13 & 4.05e-14  \\
	$ k=k_{max}$ &4.26e-14&  1.45e-16& 3.67e-14\\
	\hline 
\end{tabular}
\captionof{table}{The error between the approximated singular values $\widetilde{s}_i$ and the Exact ones $s_i$; $\left\vert \widetilde{s}_i -s_i\right\vert$ for $i=1,2,3,4$ when $m=15$.}
\label{tab last2}
\end{table}

\begin{figure}[H]
\centering
\begin{tabular}{c}
	\includegraphics[width=0.5\linewidth]{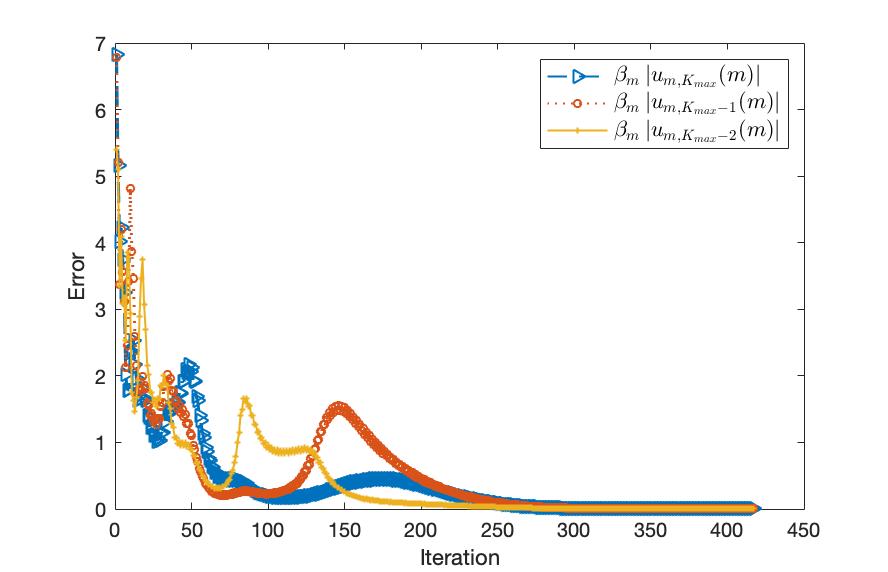}
\end{tabular}
\captionof{figure}{The error norm used to accept the last singular triplets when $k=3$ and $m=40$ on the tensor of size $50\times 20 \times 50\times 20$.}\label{fig last}
\end{figure}

In  Table, \ref{tab last3}, the Gres.norm, the cpu-time, and the number of Iterations are illustrated, when the last four singular triplets are approximated, for $m=20 \mbox{ and }30$. Figure \ref{fig last} depicts the error used to accept the three last approximated singular triplets versus the number of iterations when $m=15$ on a  $50\times 20\times  50\times 40$ tensor.
\begin{table}[H]
\centering
\small\addtolength{\tabcolsep}{-3pt}
\begin{tabular}{l|c|c|c|c|c|c|c|c|l}
	\hline size($\mathcal{A}$)&\multicolumn{3}{c|}{$100\times 100\times 50$} & \multicolumn{3}{c|}{$20\times 10\times 20\times 10$}  & \multicolumn{3}{|c}{$50\times 20\times 50\times 20$}    \\
	\hline  &GRes.norm & cpu-time & Iter &GRes.norm & cpu-time & Iter &GRes.norm & cpu-time & Iter \\
	$m=20$ &4.17e-13&6.74&3 & 2.63e-13  &10.87&102  &3.91e-12&259.59&1000   \\
	$m=30$ &3.74e-13&6.62 &2& 2.02e-13&7.13&37& 8.88e-12 &486.15 & 930   \\
	\hline
\end{tabular}
\captionof{table}{The GRes.norm, the cpu-time, and the number of Iterations (Iter) for $k=4$ and $m=20,30$.}
\label{tab last3}
\end{table}

\subsection{Tensor compression}
\begin{figure}[H]
\centering
\begin{tabular}{cccc}
	\includegraphics[width=0.2\linewidth]{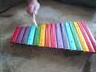}&
	\includegraphics[width=0.2\linewidth]{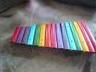}&
	\includegraphics[width=0.2\linewidth]{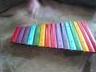}&
	\includegraphics[width=0.2\linewidth]{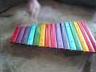}\\
\end{tabular}
\centering
\begin{tabular}{cccc}
	\includegraphics[width=0.2\linewidth]{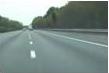}&
	\includegraphics[width=0.2\linewidth]{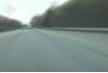}&
	\includegraphics[width=0.2\linewidth]{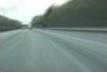}&
	\includegraphics[width=0.2\linewidth]{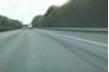}\\
	\textbf{Original} & \textbf{$k=10$} & \textbf{$k=20$} & \textbf{$k=30$}
\end{tabular}
\captionof{figure}{The $40$-th bounds of the compressed videos xylophone and departure by using Ritz method for $k=10,20,30$ with $m=40$.} \label{fig comp}
\end{figure}

In this subsection, Algorithms \ref{alg LB} and \ref{alg Ritz}
are applied to tensor compression. Algorithms \ref{alg LB} and \ref{alg Ritz} are denoted as  LB, and Ritz, respectively. The tests are carried out on color videos which are stored as fourth-order tensors. The tests are performed on the videos obtained from Matlab files : \textbf{xylophone} of size $72\times 96\times 3\times 141$ and \textbf{departure} of size $72\times 108\times 3\times 337$. To measure the quality of the compressed data, the relative error norm defined as follows is used
\[
\text{Relative error}=\dfrac{\left\Vert \mathcal{A}-\mathcal{A}_k \right\Vert_F}{\left\Vert \mathcal{A}\right\Vert_F},
\]
where $\mathcal{A}$ is the tensor representing the data, and $\mathcal{A}_k=\sum_{i=1}^{k}\widetilde{s}_k \widetilde{\mathcal{U}}_k *_2 \widetilde{\mathcal{V}}_k^T$, with $\{ \widetilde{s}_i, \widetilde{\mathcal{U}}_i, \widetilde{\mathcal{V}}_i  \}_{i=1}^k$ are the first approximated $k$ singular triplets of $\mathcal{A}$.\\
Figure \ref{fig comp} shows the $40$th bound of each the videos xylophone and departure, with the $40$th bounds of each of the compressed videos, with $k=10,20,30$, by using Ritz method with $m=40$. 

\begin{figure}[H]
\centering
\begin{tabular}{cc}
	\includegraphics[width=0.45\linewidth]{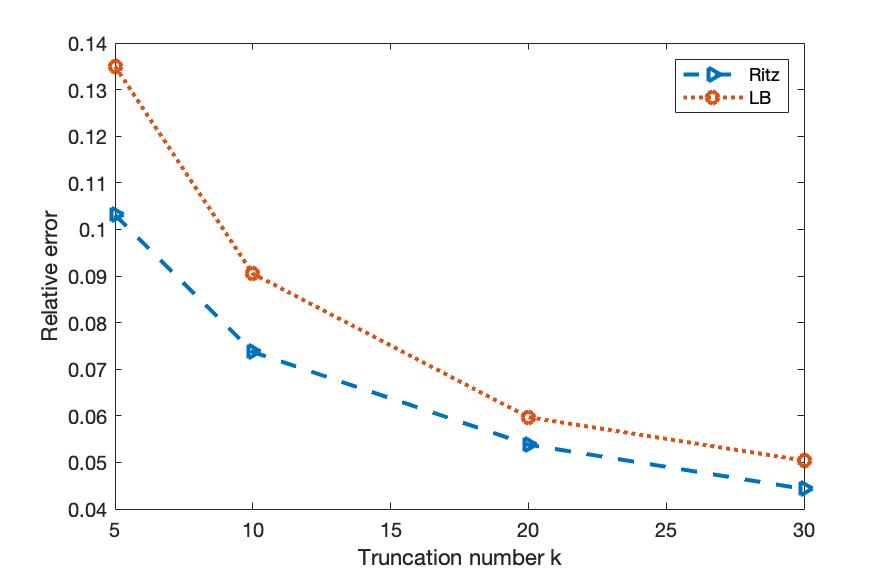}&
	\includegraphics[width=0.45\linewidth]{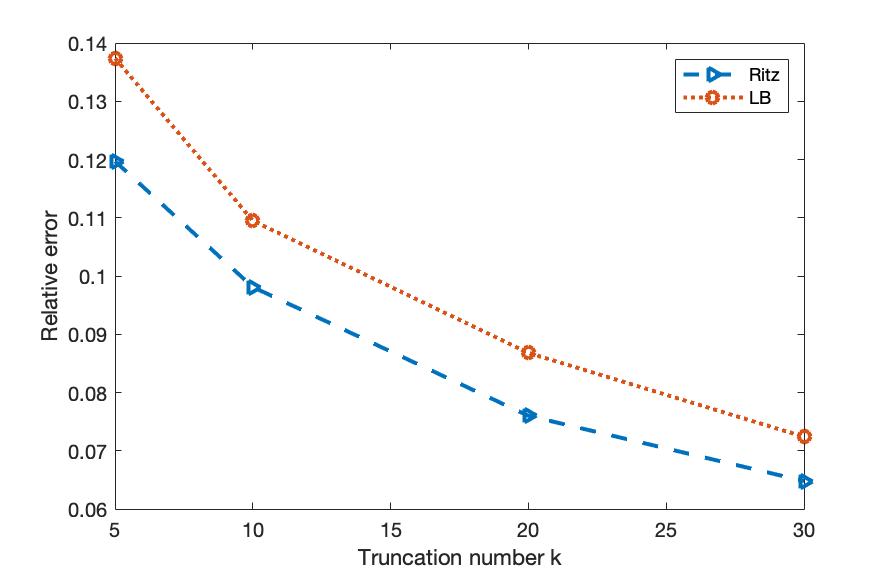}
\end{tabular}
\captionof{figure}{The relative error norm of the videos xylophone and departure, respectively from the left to the right, by using Ritz and LB methods.}\label{fig comp2}

\end{figure}

Figure \ref{fig comp2} depicts the values of the relative error norm against the values used of truncation ($k$) by using Ritz method and LB method, for both the videos xylophone and departure, respectively from the left to the right.

\noindent Figure \ref{fig comp} and \ref{fig comp2} illustrate the effectiveness of both methods for image compression, with a slight advantage to the Ritz methods.

\subsection{Facial recognition based on the Einstein product}
In this part, we illustrate some tests of facial recognition based on the Einstein product by using Algorithm \ref{alg 3}. We refer to the tests performed with Algorithm \ref{alg LB} as \textbf{LB}, and to those using Algorithm \ref{alg Ritz} as \textbf{Ritz}.\\
The used data is obtained from the Georgia Tech database GTDB crop \cite{data}. This database features the pictures of $50$ people. Each person is represented by $15$ pictures, each of which showing different facial expressions, illumination conditions, and orientation.\\
Each image is of size $50\times 50 \times 3$. These two tests are carried out:
\begin{itemize}
\item  \textbf{Test 1:} Five images of each person are randomly extracted from the database and are chosen as test images. Consequently, the database is split into two subsets : $250$  test images and $500$ training images used to form the training tensor $\widebar{\mathcal{X}}$ is of size $2500\times 3 \times 500$.
\item  \textbf{Test 2:} $3$ pictures of each person are chosen randomly as test images. The remaining $600$ images are considered as the training images. The tensor $\widebar{\mathcal{X}}$ is of size $2500\times 3 \times 600$.
\end{itemize}
The performance of both methods is measured by the identification rate (IR)
\begin{equation*}
IR=\dfrac{\text{Number of recognized images}}{\text{Number of test images}}\times 100.
\end{equation*}

\begin{figure}[H]
\centering
\begin{tabular}{ccc}
	\includegraphics[width=0.3\linewidth]{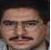}& 
	\includegraphics[width=0.3\linewidth]{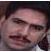}& 
	\includegraphics[width=0.3\linewidth]{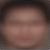}\\
	\textbf{Test image} & \textbf{Closest image} & \textbf{Mean}
\end{tabular}
\captionof{figure}{A result of \textbf{Test 1} for truncation index $k=5$.} \label{fig 2}
\end{figure}

\begin{figure}[H]
\centering
\begin{tabular}{c}
	\includegraphics[width=0.5\linewidth]{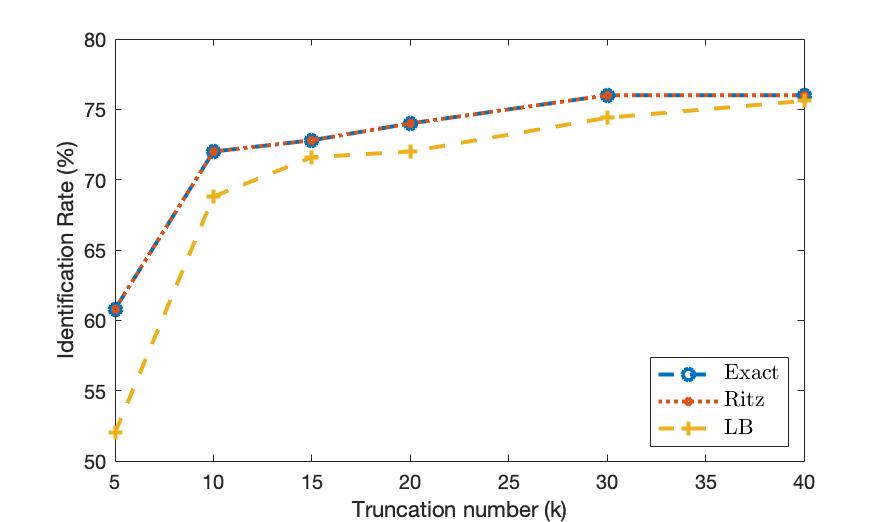}
\end{tabular}
\captionof{figure}{The curves of the identification rate (IR) versus the number of truncations used $k$, for the exact method and the approximated applied on \textbf{Test 1}.}\label{fig 3}
\end{figure}

\noindent In Figure \ref{fig 2}, we displayed one test image, the mean image and the closest image found in the database by using the Ritz method for {Test 1}, with $k=5$. In Figure \ref{fig 3} we illustrate the curves of the Identification Rate ($IR$) against the truncation index $(k)$ for  {Test 2}, by the method using the exact Einstein SVD to compute the projection space called "Exact", the LB method, and the Ritz method.

\noindent In Table \ref{tab 4}, the cpu-time is shown for the Exact, The LB, and the Ritz methods, when \textbf{Test 1} and \textbf{Test 2} are considered.
\begin{table}[H]
\centering
\begin{tabular}{c|c|c|c|c|c}
	\hline 
	Test &Method & 5 & 10 & 20 & 40 \\
	\hline \multirow{3}{*}{\textbf{Test 1}}& Exact & 295.48 & 298.93 & 307.42 & 298.50\\
	\cline{2-6} & Ritz & 57.41& 66.01 &74.97 &123.80  \\
	\cline{2-6}& LB & 49.57 &54.64 & 61.97 & 82.98\\
	\hline \multirow{2}{*}{\textbf{Test 2}}& Exact &186.30  & 208.21 & 202.02&220.93\\
	\cline{2-6} & Ritz & 42.66& 45.06 & 57.14& 97.59 \\
	\cline{2-6}& LB & 37.94 &35.77 &  40.25&52.94 \\
	\hline
\end{tabular}
\captionof{table}{The cpu-time required by the exact, the LB, and the Ritz methods for the PCA process, with different values of truncation $k=5,10,20,40$, for   {Test1} and {Test 2}.}\label{tab 4}
\end{table}
\noindent While both methods showed to perform well in terms of identification rate, we notice that the Ritz method maintains an advantage over the LB method, especially for low rank approximations, which is particularly interesting for some applications including facial recognition for which the tensor can be very large and only a few singular triplets are needed. As expected, they required less CPU-time than the Exact Einstein SVD without affecting the accuracy, with an advantage for the LB method. 

\section{Conclusion}
This paper presents an approach for approximating an extremal subset of singular triplets of a tensor using the Einstein product, regardless of its dimensionality. This approach, based on the Tensor Einstein Lanczos bidiagonalization and Ritz augmentation is particularly adapted to problems for which the computation of the complete SVD decomposition would be very time consuming. Numerical experiments show the effectiveness of this approach in terms of accuracy and CPU-time.

\end{document}